\documentclass[12pt]{article}
\usepackage{dsfont}
\usepackage{latexsym,amsfonts,amssymb}
\usepackage[dvips]{color}
\usepackage{colordvi,multicol}

\def \cal{\mathcal}

\newtheorem{thm}{Theorem}[section]

\newtheorem{pro}[thm]{Proposition}
\newtheorem{defi}{Definition}[section]
\newtheorem{rem}{Remark}[section]

\begin{document}
\title{\bf A Note on Uniform Integrability of Random Variables in a Probability Space and Sublinear Expectation Space}
\author{Ze-Chun Hu\thanks{Corresponding author.\ E-mail addresses: zchu@scu.edu.cn (Z.-C. Hu),  15266479708@163.com(Q.-Q. Zhou).}\\
 {\small College of Mathematics, Sichuan  University, Chengdu, 610064 China}\\ \\
 Qian-Qian Zhou\\
 {\small Department of Mathematics, Nanjing  University, Nanjing, 210093 China} }
\date{}
\maketitle

\begin{abstract}

In  this note we discuss uniform integrability of random variables. In a probability space, we introduce two  new notions on uniform integrability of  random variables, and prove that they are equivalent to the classic one. In a sublinear expectation space, we give de La Vall¨¦e Poussin criterion for  the uniform integrability of random variables and do some other discussions.

\end{abstract}

{\bf Key words:} uniform integrability, sublinear expectation

{\bf Mathematics Subject Classification (2000)}  60F25, 28A25

\noindent

\section{Introduction}

It is well known that the uniform integrability of a family of random variables plays an important role in probability theory. As to the uniform integrability criterions, please refer to  Chung (1974, P. 96), Chong (1979),  Chow and Teicher (1997, P. 94), Hu and Rosalsky (2011), Klenke (2014, p. 138) and Chandra (2015).

In  \cite{CHR16}, the authors introduced the notion of a sequence of random variables being {\it
uniformly nonintegrable} and gave some interesting characterizations of this uniform nonintegrability.  In \cite{HP17}, a weak notion of a sequence of random variables being uniformly nonintegrable was  introduced  and  some equivalent characterizations were given.  Motivated from \cite{CHR16} and \cite{HP17}, we will introduce two new notions of a sequence of random variables being uniformly integrable in a pobability space, and prove that they are equivalent to the classic one.

Let $(\Omega,\mathcal{F},P)$ be a probability space. Suppose that all random variables under consideration are defined on this probability space. Let $X$ be a random variable and $A\in \mathcal{F}$. We denote $E(XI_A)$ by $E(X: A)$.

\begin{defi}\label{UI}
A sequence of random variables $\{X_n,n\geq 1\}$ is said to be {\it uniformly integrable} (UI for short) if
\begin{eqnarray}\label{1}
\lim_{a\to\infty}\sup_{n\geq 1}E(|X_n|: |X_n|\geq a)=0.
\end{eqnarray}
\end{defi}

\begin{defi}\label{UNI} (\cite{CHR16})
A sequence of random variables $\{X_n,n\geq 1\}$ is said to be {\it uniformly nonintegrable} (UNI for short) if
\begin{eqnarray*}\label{1}
\lim_{a\to\infty}\inf_{n\geq 1}E(|X_n|: |X_n|\leq a)=\infty.
\end{eqnarray*}
\end{defi}

\begin{defi}\label{W-UNI} (\cite{HP17})
A sequence of random variables $\{X_n,n\geq 1\}$ is said to be {\it W-uniformly nonintegrable} (W-UNI for short) if
\begin{eqnarray*}\label{1}
\lim_{a\to\infty}\inf_{n\geq 1}E(|X_n|\wedge a)=\infty.
\end{eqnarray*}
\end{defi}

\begin{defi}\label{W*-UNI} (\cite{HP17})
A sequence of random variables $\{X_n,n\geq 1\}$ is said to be {\it W*-uniformly nonintegrable} (W*-UNI for short) if
\begin{eqnarray*}\label{defi1.1-a}
\lim_{m\to\infty}\inf_{k\geq 1}\sum_{n=0}^mP(|X_k|>n)=\infty.
\end{eqnarray*}
\end{defi}

For any random variable $X$, by the  monotone convergence theorem, we have
\begin{eqnarray*}\label{3}
\lim_{a\to \infty}E(|X|\wedge  a)=E(|X|).
\end{eqnarray*}
It follows that if $X$ is integrable, then
\begin{eqnarray}\label{3}
\lim_{a\to \infty}[E(|X|)-E(|X|\wedge  a)]=0,\ i.e. \lim_{a\to \infty}E(|X|-a: |X|\geq a)=0.
\end{eqnarray}
In virture of  (\ref{3}),  we introduce the following notion.

\begin{defi}\label{W-UI}
A sequence of random variables $\{X_n,n\geq 1\}$ is said to be {\it  W-uniformly integrable} (W-UI for short) if
\begin{eqnarray}\label{W-UI-a}
\lim_{a\to\infty}\sup_{n\geq 1}E(|X_n|-a: |X_n|\geq a)=0.
\end{eqnarray}
or equivalently,
$$
\inf_{N\geq 1}\sup_{n\geq 1}E(|X_n|-a: |X_n|\geq N)=0.
$$
\end{defi}

For any random variable $X$, we have
\begin{eqnarray}\label{1.4}
\sum_{n=1}^{\infty}P(|X|> n)\leq E|X|\leq 1+\sum_{n=1}^{\infty}P(|X|> n).
\end{eqnarray}
In virtue of (\ref{1.4}), we introduce the following notion.

\begin{defi}\label{W*-UI}
A sequence of random variables $\{X_n,n\geq 1\}$ is said to be {\it W*-uniformly integrable} (W*-UI for short) if
\begin{eqnarray}\label{W*-UI-a}
\lim_{m\to\infty}\sup_{k\geq 1}\sum_{n=m}^{\infty}P(|X_k|>n)=0.
\end{eqnarray}
\end{defi}

\begin{rem}
Let $\{X_n,n\geq 1\}$ be a sequence of random variables. It is easy to know that it is UI if and only if
 \begin{eqnarray*}\label{1}
\lim_{a\to\infty}\sup_{n\geq 1}E(|X_n|: |X_n|> a)=0.
\end{eqnarray*}
By $E(|X_n|: |X_n|> a)=E|X_n|-E(|X_n|: |X_n|\leq  a)$, we can say that UI corresponds to W-UI in some sense.
Similarly, we can say that  W-UI corresponds to W-UNI and  W*-UI corresponds to W*-UNI in some sense, respectively.
\end{rem}

In Section 2, we will prove that UI, W-UI and W*-UI are equivalent in a probability space.

Recently, motivated by the risk measures,
 superhedge pricing and modeling uncertain in finance, Peng \cite{Pe04}-\cite{Pe10}
 initiated the notion of independent and identically distributed (IID) random
 variables under sublinear expectations,  proved the weak law of large
 numbers and the central limit theorems, defined the $G$-expectations, $G$-Brownian motions and built It\^{o}'s type stochastic calculus.
 In Section 3, we discuss uniform integrability of  random variables in a sublinear expectation space, and  present de La Vall¨¦e Poussin criterion for  the uniform integrability of random variables and make some other discussions.

\section{Uniform integrability in a  probability space}\setcounter{equation}{0}

In \cite{HP17}, we prove that
$$
\mbox{UNI} \Rightarrow \mbox{W-UNI} \Leftrightarrow \mbox{W*-UNI},
$$
and W-UNI is strictly weaker than UNI in general.  While, as to UI, W-UI and W*-UI, we have the following result.

\begin{thm}\label{thm}
\begin{eqnarray}\label{thm-a}
\mbox{UI} \Leftrightarrow \mbox{W-UI} \Leftrightarrow \mbox{W*-UI}.
\end{eqnarray}
\end{thm}
{\bf Proof.}  Let $\{X_n,n\geq 1\}$ be a sequence of random variable in a probability space $(\Omega,\mathcal{F},P)$.

 UI $\Rightarrow$ W-UI: Suppose that $\{X_n,n\geq 1\}$ is UI.  Then by Definition \ref{UI}, Definition \ref{W-UI}
 and the inequality
 $$E[|X_n|-a:|X_n|>a]\leq E[|X_n|:|X_n|>a],$$
 we know that $\{X_n,n\geq 1\}$ is W-UI.

 W-UI $\Rightarrow$ UI:  Suppose that $\{X_n,n\geq 1\}$ is W-UI.
 For any set $A\in \cal{F}$,  any positive constant $C$ and any integer $n$, we have
\begin{eqnarray}\label{proof-a}
\int_A |X_n|dP&=&\int_{A\cap [|X_n|\geq C]} (|X_n|-C)dP+CP(A\cap [|X_n|\geq C])+\int_{A\cap [|X_n|< C]} |X_n|dP\nonumber\\
&\leq&\int_{[|X_n|\geq C]} (|X_n|-C)dP+2CP(A)\nonumber\\
&\leq&\sup_{k\geq 1}\int_{[|X_k|\geq C]} (|X_k|-C)dP+2CP(A).
\end{eqnarray}
By the definition of W-UI, there exists a positive number $C_0$ such that
\begin{eqnarray}\label{proof-b}
\sup_{k\geq 1}\int_{[|X_k|\geq C_0]} (|X_k|-C_0)dP<\frac{\epsilon}{2}.
\end{eqnarray}
Let $\delta=\frac{\epsilon}{4C_0}$.  Then for any $A\in \cal{F}$ with $P(A)<\delta$, by
(\ref{proof-a}) and (\ref{proof-b}), we obtain that
\begin{eqnarray}\label{proof-c}
\int_A |X_n|dP<\epsilon,\ \forall n\geq 1.
\end{eqnarray}
Setting $A=\Omega$ in (\ref{proof-a}) and using the definition of W-UI, we get that
\begin{eqnarray}\label{proof-d}
\sup_{n\geq 1}E(|X_n|)<\infty.
\end{eqnarray}
By (\ref{proof-c}) and (\ref{proof-d}), we obtain that $\{X_n,n\geq 1\}$ is UI.

W-UI $\Rightarrow$ W*-UI: For any random variable $X$ and any positive integer $m$, by Fubini's theorem, we have
\begin{eqnarray*}
\sum_{n=m}^{\infty}P(|X|>n)&=& \sum_{n=m}^{\infty}\int_{n}^{n+1}P(|X|>n)dx\\
&\leq&\sum_{n=m}^{\infty}\int_{n}^{n+1}P(|X|>x-1)dx\\
&=&\int_m^{\infty}P(|X|>x-1)dx\\
&=&\int_{m-1}^{\infty}P(|X|>x)dx\\
&=&\int_{\Omega}\left(\int_{m-1}^{\infty}I_{\{|X|>x\}}dx\right)dP\\
&=&E(|X|-(m-1): |X|>m-1)\\
&=&E(|X|-(m-1):|X|\geq m-1).
\end{eqnarray*}
It follows that W-UI $\Rightarrow$ W*-UI.

W*-UI $\Rightarrow$ W-UI: For any random variable $X$ and any positive integer $m$, by Fubini's theorem, we have
\begin{eqnarray*}
\sum_{n=m}^{\infty}P(|X|>n)&\geq & \sum_{n=m}^{\infty}\int_{n}^{n+1}P(|X|>x)dx\\
&=&\int_{m}^{\infty}P(|X|>x)dx\\
&=&\int_{m}^{\infty}\left(\int_{\Omega}I_{\{|X|>x\}}dP\right)dx\\
&=&\int_{\Omega}\left(\int_{m}^{\infty}I_{\{|X|>x\}}dx\right)dP\\
&=&E(|X|-m: |X|>m)\\
&=& E(|X|-m: |X|\geq m).
\end{eqnarray*}
It follows that W*-UI $\Rightarrow$ W-UI.

Hence (\ref{thm-a}) holds, and the proof is complete. \hfill\fbox

\section{Uniform integrability in a  sublinear expectation space}\setcounter{equation}{0}

In this section, we discuss the uniform integrability of random variables in a  sublinear expectation space. At first,  we present some basic settings about sublinear expectations. Please refer to  Peng \cite{Pe04}-\cite{Pe10}, and Cohen et al. \cite{CJP11}  for more details.

Let $(\Omega,\cal{F})$ be a given measurable space and  $\bf{\mathcal {H}}$ be a linear space of $\cal{F}$-measurable real functions defined on $\Omega$ such that  for any constant number $c, c\in \mathcal{H}$;  if $X\in\mathcal{H}$, then $|X|\in\mathcal{H}$ and $XI_A\in \cal{H}$ for any $A\in \cal{F}$.

\begin{defi}
A  sublinear expectation $\cal{E}$ on $\mathcal{H}$ is a functional $\cal{E}:\mathcal{H}\rightarrow\mathbb{R}$ satisfying the
following properties:\\
\hspace*{0.3cm} (a) Monotonicity: $\cal{E}[X]\ge\cal{E}[Y]$, if $X\ge Y$.\\
\hspace*{0.3cm}  (b) Constant preserving: $\cal{E}[c]=c,\forall c\in\mathbb{R}.$\\
\hspace*{0.3cm}  (c) Sub-additivity: $\cal{E}[X+Y]\le\cal{E}[X]+\cal{E}[Y].$\\
\hspace*{0.3cm} (d) Positive homogeneity: $\cal{E}[\lambda X]=\lambda\cal{E}[X]$, $\forall\lambda\ge 0.$\\
The triple $(\Omega,\mathcal{H},\cal{E})$ is called a
sublinear expectation space.
\end{defi}

\begin{defi}(\cite[Definition 3.1]{CJP11})
For $p \in [1, \infty)$, the map
$$
{|| \cdot ||}_p : X \mapsto (\cal{E}[|X|^p])^{1/p}
$$
forms a seminorm on $\cal{H}$. Define the space ${\cal{L}}^{p}(\cal{F})$ as the completion under ${||\cdot||}_{p}$ of the set
$$
\{X\in \cal{H}: {||X||}_{p} < \infty\}
$$
and then $L^{p}(\cal{F})$ as the equivalence classes of ${\cal{L}}^{p}$ modulo equality in ${||\cdot||}_{p}$.
\end{defi}

\begin{defi}\label{defi2.4}(\cite[Definition 3.2]{CJP11})
Consider $K\subset L^{1}$. $K$ is said to be uniformly integrable
if $\lim_{c\to\infty}\sup_{X\in K}\cal{E}[I_{\{|X|\geq c\}}|X|]=0$.
\end{defi}


\begin{thm}\label{thm3.1}(\cite[Theorem  3.1]{CJP11})
Suppose $K$ is a subset of $L^1$. Then $K$ is uniformly integrable if
and only if the following two conditions hold.
\begin{itemize}
\item[(i)]  $\{\cal{E}[|X|]\}_{X\in K}$ is bounded.

\item[(ii)]  For any $\epsilon<0$  there is a $\delta>0$ such that for all $A \in \cal{F}$ with $\cal{E}[I_A]\leq \delta$, we
have $\cal{E}[I_A|X|] < \epsilon$  for all $X\in K.$
\end{itemize}
\end{thm}

Now we present the following de La Vall¨¦e Poussin criterion for  the uniform integrability.

\begin{thm}\label{thm3.2}
Let $K$ be a subset of $L^1$. Then $K$ is uniformly integrable if and only if  there is a nonnegative function
$\varphi$ defined on $[0,\infty)$ such that $\lim_{t\to\infty}\varphi(t)/t = \infty$ and $\sup_{X\in K} \cal{E}[\varphi\circ |X|]< \infty$.
\end{thm}
{\bf Proof.} As to the sufficiency, refer to \cite[Corollary 3.1.1]{CJP11}.  In the following, we give the proof of the necessity. The idea comes from the corresponding proof in a probability space (see e.g. \cite[Theorem 7.4.5]{Yan}).

Suppose that $K$ is uniformly integrable. For any constant $a>0$, we have $\cal{E}[(|X|-a)^+]\leq \cal{E}[|X|I_{\{|X|\geq a\}}]$.  It follows that there exists a sequence $\{n_k\}$ of integers such that $n_k\uparrow \infty$  and
\begin{eqnarray}\label{thm3.2-a}
\sup_{X\in K}\cal{E}[(|X|-n_k)^+]<2^{-k},\ k\geq 1.
\end{eqnarray}
Define a function
$$
\varphi(t)=\sum_{k\geq 1}(n-n_k)^+,\ n\leq t<n+1,\ n=0,1,2,\cdots.
$$
Then $\varphi$ is a nonnegative, nondecreasing and right continuous function. What's more, we have
$$
\lim_{n\to\infty}\frac{\varphi(n)}{n}=\lim_{n\to\infty}\sum_{k\geq 1}(1-\frac{n_k}{n})^+=\infty,
$$
which implies that $\lim_{t\to\infty}\varphi(t)/t = \infty$.

By Fubini's theorem, the monotone convergence theorem (\cite[Theorem 2.2]{CJP11}),  the sublinear property of $\cal{E}$ and (\ref{thm3.2-a}), we obtain that for any $X\in K$,
\begin{eqnarray*}
\cal{E}[\varphi\circ |X|]&=&\cal{E}\left[\sum_{n=0}^{\infty}\sum_{k=1}^{\infty}(n-n_k)^+I_{\{n\leq |X|<n+1\}}\right]\\
&=&\cal{E}\left[\sum_{k=1}^{\infty}\sum_{n=0}^{\infty}(n-n_k)^+I_{\{n\leq |X|<n+1\}}\right]\\
&=&\cal{E}\left[\lim_{m\to\infty}\sum_{k=1}^m\sum_{n=0}^{\infty}(n-n_k)^+I_{\{n\leq |X|<n+1\}}\right]\\
&=&\lim_{m\to\infty}\cal{E}\left[\sum_{k=1}^m\sum_{n=0}^{\infty}(n-n_k)^+I_{\{n\leq |X|<n+1\}}\right]\\
&\leq&\lim_{m\to\infty}\sum_{k=1}^m\cal{E}\left[\sum_{n=0}^{\infty}(n-n_k)^+I_{\{n\leq |X|<n+1\}}\right]\\
&=&\sum_{k=1}^{\infty}\cal{E}[(|X|-n_k)^+]<1.
\end{eqnarray*}
 \hfill\fbox\\

With respect to Definitions 1.5 and 1.6, we introduce the following two notions.

\begin{defi} \label{W-UI-SLE}
Consider $K\subset L^1.$ $K$ is said to be W-uniformly integrable (W-UI for short)
if
\begin{eqnarray}\label{W-UI-3-a}
\lim_{a\to\infty}\sup_{X\in K}\cal{E}[I_{\{|X|\geq a\}}(|X|-a)]=0.
\end{eqnarray}
\end{defi}

\begin{defi} \label{W*-UI-SLE}
Consider $K\subset L^1.$ $K$ is said to be S-uniformly integrable (S-UI for short)
if
\begin{eqnarray}\label{W*-UI-3-a}
\lim_{m\to\infty}\sup_{X\in K}\sum_{n=m}^{\infty}\cal{E}[I_{\{|X|>n\}}]=0.
\end{eqnarray}
\end{defi}

\begin{pro}\label{pro-3.3}
Suppose that $K$ is a family of random variables in a sublinear expectation $(\Omega,\cal{H},\cal{E})$. Then we have
\begin{eqnarray}\label{pro-3.3-a}
\mbox{UI} \Leftrightarrow \mbox{W-UI} \Leftarrow \mbox{S-UI}.
\end{eqnarray}
\end{pro}
{\bf Proof.} UI $\Rightarrow$ W-UI:  Suppose that $K$ is UI.  Then by Definition \ref{defi2.4}, Definition \ref{W-UI-SLE}
 and the inequality
 $$\cal{E}[(|X|-a)I_{\{|X|>a\}}]\leq E[|X|I_{\{|X|>a\}}],$$
 we know that $K$ is W-UI.

 W-UI $\Rightarrow$ UI:  Suppose that $K$ is W-UI.
 For any set $A\in \cal{F}$,  any positive constant $C$ and any  $X\in K$, we have
\begin{eqnarray}\label{pro-3.3-b}
\cal{E}[|X|I_A]&=&\cal{E}\left[(|X|-C)I_{A\cap \{|X|\geq C\}}+CI_{A\cap \{|X|\geq C\}}+|X|I_{A\cap \{|X|<C\}}\right]\nonumber\\
&\leq&\cal{E}[(|X|-C)I_{A\cap \{|X|\geq C\}}]+\cal{E}[CI_{A\cap \{|X|\geq C\}}]+\cal{E}[|X|I_{A\cap \{|X|<C\}}]\nonumber\\
&\leq&\cal{E}[(|X|-C)I_{\{|X|\geq C\}}]+2C\cal{E}[I_A].
\end{eqnarray}
By Definition \ref{W-UI-SLE} there exists a positive number $C_0$ such that
\begin{eqnarray}\label{pro-3.3-c}
\sup_{X\in K}\cal{E}[(|X|-C)I_{\{|X|\geq C\}}]<\frac{\epsilon}{2}.
\end{eqnarray}
Let $\delta=\frac{\epsilon}{4C_0}$.  Then for any $A\in \cal{F}$ with $\cal{E}[I_A]<\delta$, by
(\ref{pro-3.3-b}) and (\ref{pro-3.3-c}), we obtain that
\begin{eqnarray}\label{pro-3.3-d}
\cal{E}[|X|I_A]<\epsilon,\ \forall X\in K.
\end{eqnarray}
Setting $A=\Omega$ in (\ref{pro-3.3-b}) and using Definition \ref{W-UI-SLE}, we get that
\begin{eqnarray}\label{pro-3.3-e}
\sup_{X\in K}\cal{E}[|X|]<\infty.
\end{eqnarray}
By (\ref{pro-3.3-d}), (\ref{pro-3.3-e}) and Theorem \ref{thm3.1}, we obtain that $K$ is UI.

S-UI $\Rightarrow$ W-UI: Suppose that $K$ is S-UI. For any $X\in K$ and any integer $m$, by the monotone convergence theorem (\cite[Theorem 2.2]{CJP11}) and  the sublinear property of $\cal{E}$, we get
\begin{eqnarray*}
\cal{E}[(|X|-m)I_{\{|X|\geq m\}}]&=&\cal{E}[(|X|-m)I_{\{|X|> m\}}]\\
&=&\cal{E}\left[\int_m^{\infty}I_{\{|X|> x\}}dx\right]\\
&=&\cal{E}\left[\sum_{n=m}^{\infty}\int_n^{n+1}I_{\{|X|> x\}}dx\right]\\
&=&\cal{E}\left[\lim_{l\to\infty}\sum_{n=m}^l\int_n^{n+1}I_{\{|X|> x\}}dx\right]\\
&=&\lim_{l\to\infty}\cal{E}\left[\sum_{n=m}^l\int_n^{n+1}I_{\{|X|> x\}}dx\right]\\
&\leq&\lim_{l\to\infty}\sum_{n=m}^l\cal{E}\left[\int_n^{n+1}I_{\{|X|> x\}}dx\right]\\
&\leq&\sum_{n=m}^{\infty}\cal{E}\left[I_{\{|X|> n\}}\right],
\end{eqnarray*}
which together with Definitions \ref{W-UI-SLE} and \ref{W*-UI-SLE} implies that $K$ is W-UI. \hfill\fbox

\begin{rem}
The part ``W-UI $\Rightarrow$ W*-UI"of the proof of Theorem \ref{thm} tell us that in general we don't have that W-UI $\Rightarrow$ S-UI in a sublinear expecation space. In the following, we will give a counterexample.

Let $\Omega=\{0,1,2,\ldots\}$. For any $n=2,\ldots,$ define a probability measure $P_n$ on $\Omega$ as follows:
$$
P_n(n)=\frac{1}{n\ln n},\ \  P_n(0)=1-\frac{1}{n\ln n}.
$$
Denote by $E_n$ the expectation with respect to the probability measure $P_n$. Define the sublinear expectation $\cal{E}$ by
$$
\cal{E}[\cdot]:=\sup_{n\geq 2}E_n[\cdot].
$$
Let $X$ be a random variable defined on $\Omega$ by
$$
X(n)=n,\ n=0,1,2,\ldots.
$$

We have
\begin{eqnarray*}
\lim_{m\to\infty}\cal{E}[|X|I_{\{|X|\geq m\}}]&=&
\lim_{m\to\infty}\sup_{n\geq 2}E_n[|X|I_{\{|X|\geq m\}}]\nonumber\\
&=&\lim_{m\to\infty}\sup_{n\geq m}n\times \frac{1}{n\ln n}\nonumber\\
&=&\lim_{m\to\infty}\frac{1}{\ln m}=0,
\end{eqnarray*}
which implies that $\{X\}$ is UI and thus $\{X\}$ is W-UI by Proposition 3.3.

We also have
\begin{eqnarray*}
\lim_{m\to\infty}\sum_{n=m}^{\infty}\cal{E}[I_{\{|X|>n\}}]&=&
\lim_{m\to\infty}\sum_{n=m}^{\infty}\sup_{k\geq 2}E_k[I_{\{|X|>n\}}]\\
&=&\lim_{m\to\infty}\sum_{n=m}^{\infty}\sup_{k\geq n+1}\frac{1}{k\ln k}\\
&=&\lim_{m\to\infty}\sum_{n=m}^{\infty}\frac{1}{(n+1)\ln (n+1)}=+\infty,
\end{eqnarray*}
which implies that $\{X\}$ is not S-UI.
\end{rem}

\bigskip

{ \noindent {\bf\large Acknowledgments} \quad The authors  thank the anonymous referee for providing helpful comments to improve  the manuscript.  This work was supported by National Natural Science Foundation of China (Grant No. 11371191).}

\end{document}